
\input vanilla.sty
\magnification 1200
\baselineskip 18pt
\overfullrule=0pt
\input definiti.tex
\input mathchar.tex
\define\pmf{\par\medpagebreak\flushpar}

\define\pbf{\par\bigpagebreak\flushpar}

\pmf
\title
Fundamental Group of Self-Dual Four-Manifolds with Positive \\
Scalar Curvature
\endtitle

\author
Alexander G. Reznikov
\endauthor
\author
June, 1994
\endauthor
\author
Revised March, 1995
\endauthor

\subheading{1. Background}  It is a fundamental goal in Riemannian geometry
to understand the topology of manifolds of positive curvature.  The
only general facts so far known are: finiteness of fundamental group (Myers'
theorem), vanishing of $\hat A$-genus (Lichnerowicz' theorem and the
modification of Hitchin) and a universal bound on Betti numbers (Gromov's
theorem).  In a well-known paper [10] Micallef and Moore introduced
a new notion of positivity for the curvature tensor, that is,
{\it positivity on complex isotropic two-planes}.  For $x \in M$, a
Riemannian manifold,
let $R : \Lambda^2 T_x M \to \Lambda^2 T_x M$ be the curvature tensor.
After complexification we get a Hermitian operator in
$\Lambda^2 T^\bbc_x M$.  We say that $z \in \Lambda^2 T^\bbc_x M$ comes
from a complex isotropic
two-plane if $z = \xi \wedge \eta$ with $(\xi, \xi)_\bbc =
(\xi, \eta)_\bbc = (\eta, \eta)_\bbc = 0$.  Here $(\cdot, \cdot)_\bbc$ is the
canonical symmetric (not Hermitian!) complexification of the Euclidean
scalar product in $T_x M$.  The condition above says that $(R z, z) > 0$
for such
$z$.  The theorem of Micallef and Moore reads:

\proclaim{\bf 1.1. Theorem}  (Micallef, Moore).  Let $M^n$ be a compact, simply
connected Riemannian manifold.  If $(R z, z) > 0$ for $z$ as above, then $M$
is homeomorphic to $S^n$.
\endproclaim

A remarkable application is a pointwise pinching theorem, which reads as
follows:

\proclaim{\bf 1.2. Corollary}  Let $M$ be a compact Riemannian manifold, whose
sectional curvature satisfies $\frac14 B (x) < K (x) \le B (x)$ for some
positive function $B$, then $M \approx S^n$.
\endproclaim

If one allows equality in $(R z, z) \ge 0$ or $\frac14 B (x) \le K (x) \le B
(x)$, one
encounters more topological types of manifolds, like complex
projective spaces.

On the other hand LeBrun [2b] classified the underlying topological manifolds
of simply-connected four-manifolds of positive scalar curvature.

In both cases, given the  absence of the  Myers theorem, cited above, one asks
an
important question:
\pmf
{\bf 1.4 Which fundamental group may manifolds with such a curvature
have?}

Below in Theorem 3.3 we will derive a very strong restriction on $\pi_1 (M)$ in
case
of dimension four, namely:

\proclaim{\bf Main Theorem (3.3)}  Let $M$ be a compact four-dimensional
manifold
either with curvature, positive on complex isotropic two-planes, or self-dual
of positive scalar curvature.  If
$\pi_1 (M)$
admits a nontrivial unitary representation, and $M$ is orientable, then there
exists
a surjective homomorphism from $\pi_1 (M)$ on $\bbz$.
\endproclaim

\proclaim{\bf Corollary}  If $\pi_1 (M)$ is finite, then either $\pi_1 (M) =
1$,
or $\pi_1 (M) = \bbz_2$.
\endproclaim

Observe that finitely presented groups which do not admit
a nontrivial unitary representation, are extremely rare (see 3.4).

We will also discuss the diffeomorphism problem, namely
\pmf
{\bf 1.5.  Is it true,
that in conditions of Theorem 1.1, $M$ is diffeomorphic to $S^n$?}

In connection to this question, recall that some exotic spheres admit
a metric of positive sectional curvature [8].  On the other
hand, ``strongly'' pinched manifolds are standard spheres [9] , [11], [15].
Here, we suggest a completely new approach to this problem, using the theory of
SD-connections
on four-manifolds [6], especially the Donaldson's collar theorem.

\subheading{2. Positivity of curvature on complex isotropic two-planes: an
algebraic study}
\pmf
2.1  The computation that follows is essentially contained in the Micallef and
Moore's  paper.Let $V$ be a four-dimensional orientable Euclidean space.
Consider
the canonical decomposition $\Lambda^2 V = \Lambda^2_+  V \oplus
\Lambda^2_ - V$. The
unit sphere $S^2_{\pm} (V)$ in $\Lambda^2_{\pm} (V)$ consists of
twistors, that
is, orthogonal complex structures in $V$.  Let $Z \subset V^{\bbc}$ be a
a complex
isotropic two-plane.  Since $Z \cap V = Z \cap i V = 0$, we may
look at $Z$ as a graph of
an invertible real operator $\Cal P: V \to V$.  The equation
$(v + i \Cal P v, w + i \Cal P w) = 0$ for any
$v, w \in V$ implies $\Cal P \in S^2_{\pm} (V)$.
Next, let $R: \Lambda^2 V \to \Lambda^2 V$ be an
curvature-like symmetric operator.  We denote the Hermitian
complexification of $R$ again by $R$.  Choose unit vectors
$v, w$ in $V$ such that
$(v, w) = (\Cal P v, w) = 0$.  The complex plane, corresponding to
$\Cal P$ is spanned
by $v + i \Cal P v$ and $w + i \Cal P w$, so the element $z$ in
$\Lambda^2_{\bbc} V$ is
$(v \wedge w - \Cal P v \wedge \Cal P w) + i (\Cal P v \wedge w + v \wedge
\Cal P w)$.
Assume $\Cal P \in S^2 _+ (V)$.  Then the bivectors $f = \frac1{\sqrt{2}}
(v \wedge w - \Cal P v \wedge \Cal P w)$ and $g = \frac1{\sqrt{2}}
(\Cal P v \wedge
w + v \wedge \Cal P w)$ are
complementing elements of an orthonormal basis $(f, g, \Cal P)$ of $\Lambda^2_+
(V)$.
The condition $(R z, z) \ge 0$ therefore reads $(R f, f) + (R g, g)
\ge 0$, or $Tr R |_{\Lambda^2_+ (V)} \ge R (\Cal P, \Cal P)$.
We will state this in a form of lemma, in which
$s = 4 T r R|_{\Lambda^2_{\pm} (V)}$ is scalar
curvature, and $W$ is the Weyl tensor.

\proclaim{\bf 2.2 Lemma}  A curvature-like tensor $R: \Lambda^2 V \to \Lambda^2
V$ is nonnegative
on complex isotropic two-planes, if and only if for a unit vector
$\Cal P$ of $\Lambda^2_{\pm} (V)$ one has $(W \Cal P, \Cal P) \le \frac{s}{6}$.
\endproclaim

\subheading{3. Vanishing results and the fundamental group.}
\pmf
{\bf 3.1}  Consider the modified de Rham complex [6],
$$ \Omega^0 (M) \overset d \to \rightarrow \Omega^1 (M) \overset \Cal D \to
\rightarrow \Omega^2_{\pm} (M),$$
where $\Omega^i (M)$ stands for the sheaf of $C^\infty i$-forms on $M$, and
$\Omega^2_{\pm} (M)$ is a sheaf of sections of $\Lambda^2_{\pm} T^\ast (M)$.

The Weitzenb\"ock formula [1], [6] gives
$$ \Cal D^\ast \Cal D = \nabla^\ast \nabla - 2 W^{\pm} + \frac{s}{3} $$
Accounting 2.2, we conclude (see also [14]).

\proclaim{\bf Proposition  (3.1)}  Let $M$ be a compact four-manifold with a
curvature, non-negative on complex
isotropic two-planes, or a self-dual manifold with non-negative  scalar
curvature.  Then any self-dual or anti-selfdual harmonic two-form on $M$ (resp.
anti-selfdual form) is self-parallel.  If
the curvature above is positive, then $H^2 (M, \bbr) = 0$.(resp. $H^2_-(M,
\bbr)=0$).
\endproclaim
\pmf
{\bf 3.2}  Here we will derive twisted versions of 3.1.  Let $\rho: \pi_1 (M)
\to U (n)$ be a unitary
representation.  Let $E_\rho$ be the corresponding flat Hermitian vector bundle
over $M$.  Let $(\Omega^\ast (M, E), d)$ be the complex
of $E$-valued forms, where $d$ is induced by the flat connection.  The
cohomology of $\Omega^\ast (M, E)$ coincides with the cohomology of
the local system associated to $\rho$.  One has the Laplace operator
$\Delta_i$ aching in $\Omega^i (M, E)$
and the Hodge theorem\break Ker $\Delta_i \approx H^i (M, E)$.  The proof of
the following
theorem is parallel to the proof of 3.1.

\proclaim{\bf Theorem (3.2)}  Let $M$ be as in 3.1.  Then any self-dual or
antiself-dual
harmonic two-form with coefficients in $E$ (resp. antiself-dual
harmonic two-form with coefficients in $E$)  is self-parallel.  If the
curvature on complex isotropic two-planes is positive (resp. the scalar
curvature is positive), then $H^2 (M, E) = 0$ (resp. $H^2_- (M, E) = 0$).
\endproclaim
\pmf
{\bf 3.3.}  Now, assume that $M$ is compact four-manifold with curvature,
positive on
complex isotropic two-planes or a selfdual manifold with positive scalar
curvature.  Then by Theorem 3.1, $H^2 (M, \bbr) = 0$ ( resp.$H^2_- (M, \bbr) =
0$)  and by
Theorem 3.2., $H^2 (M, E_\rho) = 0$ ( resp.  $H^2_- (M, E_\rho) = 0$) for any
unitary representation $\rho$.  This
implies the following strong statement concerning the structure of $\pi_1 (M)$.

\proclaim{\bf Main Theorem (3.3)}  Suppose $M$ is orientable and there exists a
nontrivial
unitary representation of $\pi_1 (M)$.  Then there exists a surjective
homomorphism
$\pi_1 (M) \to \bbz$.
\endproclaim

\demo{\bf Proof} Suppose $H^1 (M, \bbr) = 0$ and let $\rho: \pi_1 (M) \to U
(n)$ be
nontrivial.  Then $h^0 (M, E_\rho) < n$.  The Poincar\'e duality and Theorem
3.2. give $\chi (M, E) = 2 h^0 (M, E_\rho) - 2 h^1 (m, E_\rho)$.  Similarly,
$\chi (M) = 2 - 2 b_1 (M)$.  Since $b_1 (M) = 0$, we get
$\chi (M, E) < n \chi (M)$.  On the
other hand, $c_i (E_\rho) = 0$ because $E_\rho$ is flat, so by $A S$ index
theorem we get $\chi (M, E) = n \chi (M)$.  This contradiction shows that
$H^1 (M, \bbr) \neq  0$, hence the result.The proof for the selfdual manifolds
with positive scalar curvature is parallel, but one looks at the modified
complex
$\Omega^o (E_\rho) \overset d \to \rightarrow
\Omega^1 (E_\rho) \overset \Cal D \to \rightarrow \Omega^2_- (E_\rho)$.

\proclaim{\bf 3.4 Corollary}  In conditions of 3.3, let $\pi_1 (M)$ be finite.
Then
either $\pi_1 (M) = 1$, or $\pi_1 (M) = \bbz_2$.
\endproclaim

\demo{\bf Proof}  Consider the orientable covering of $M$ and apply Theorem
3.3.

The f.p. groups which do not yield conditions of Theorem 3.3 are rare.  They
may not admit any nontrivial linear representation over any field
and any subgroup of finite index.

\subheading{4. Concluding remarks}
\pmf
{\bf 4.1.}  Here we address the problem, stated in 1.4:

Is any compact simply-connected four-manifold, with curvature, positive
on complex isotropic two-planes, diffeomorphic to a four-sphere?

A weaker question is:
\pmf
{\bf 4.2.}  Is any compact four-manifold $M$, with pointwise $\frac14$-pinched
curvature
diffeomorphic to $S^4$?

There are several ways to prove that a $\delta$-pinched manifold
with $\delta$ sufficiently closed to 1, is a standard sphere.  [9], [11].
We sketch here an approach to settle 1.4, as follows.
\pmf
{\bf 4.3.}  Fix a $S U (2)$-principal bundle $\eta$ over $M$ with $(c_2, [M]) =
-1$.  Consider
the moduli space $\Cal M$ of $S D$-connections in $\eta$.  Then the
following result holds.

\proclaim{\bf Lemma (5.3)}  The metric of $M$ behaves like a generic
metric, that is, the canonical compactification
$\bar {\Cal M}$ of $\Cal M$ is a
smooth compact 5-manifold with $M$ as a boundary.
\endproclaim

\demo{\bf Proof}  The standard analysis shows ([6], p.69 ) that for $\Cal M$
to be a smooth manifold it is enough that for any $S D$-connection $\nabla$,
the
operator
$$ \Cal D: \Omega^1 (ad \, \eta) \to \Omega^2_- (ad \, \eta) $$
will be surjective.  Since $\Cal D {\Cal D}^\ast = \nabla^\ast \nabla -
2 W^- + s/3$ for $\nabla$ self-dual
([6], p.111), and since $s/3- 2 W^-$ is a positive operator
by lemma $\Cal M$ is smooth.  Since $H^2 (M, \bbr) = 0$ there are no
reducible connections, so $\partial \bar {\Cal M} = M$ by the Donaldson's
theorem.
\pmf
{\bf 4.4.}  Now, there is a canonical Finsler metric on $\Cal M$, coming from
$L^4$-metric on $\Omega^1 (ad \, \eta)$.  This metric is invariant under
conformal
changes of metric on $M$.  In case $M = (S^4, can)$ this Finsler metric
is in fact a (Riemannian) hyperbolic metric of $B^5$.  In any case,
by standard methods of Finsler manifolds, one finds a canonical
``osculating'' Riemannian metric $g$ in $\Cal M$.  The direct
computation shows the following:

\proclaim{\bf Lemma (4.4)}  If $M$ is $\delta$-pinched with $\delta$ close to
1, then
the curvature of $({\Cal M}, g)$ is pinched negative.
\endproclaim

Then it follows immediately that $M = \partial \bar {\Cal M} \approx S^4$.
We hope that a finer analysis will be helpful to settle 1.4 and 5.1.

\pbf
\centerline{\bf References}

\item{1.} J. - P. Bourguignon, Formules de Weitzenb\"ock en dimension 4,
in:
Arthur Besse, G\'eometrie Riemannienne en Dimension 4, Paris, Cedic/ Fernand
Nathan, 1981.

\item{2a.} C. LeBrun, Explicit self-dual metrics on $\bbc P^2
\# \cdots \# \bbc P^2$, J. Diff. Geom., {\bf 34} (1991),
223--254.

\item{2b.}C. LeBrun, On the topology of self-dual 4-manifolds, Proc. Amer.
Math.Soc., {\bf 98} (1986), 637--642.

\item{3.} S. K. Donaldson, An application of gauge theory to four dimensional
topology, J.
Diff. Geom., {\bf 18} (1983), 279--315.

\item{4.} S. K. Donaldson, R. Friedman, Connected sums at self-dual
manifolds
and deformations of singular spaces, Nonlinearity, {\bf 2} (1989), 197--239.

\item{5.} A. Floer, Self-dual conformal structures on
$\ell \bbc P^2$, J. Diff. Geom., {\bf 33} (1991), 551--574.

\item{6.} D. Freed, K. Uhlenbeck, Instantons and Four-Manifolds,
Springer-Verlag, 1984.

\item{7.} M. Friedman, The topology of four-dimensional manifolds, J. Diff.
Geom.,
{\bf 17} (1982), 357--454.

\item{8.} D. Gromoll, W.Meyer, An exotic sphere with nonnegative sectional
curvature,
Annals of
Math.{\bf 100} (1974), 401--406

\item{9.} C. Margerin, Pointwise pinched manifolds are space forms,
Proc. Symp. Pure
Math. {\bf 44} (1986), 307--328.

\item{10.} M. Micallef, J. D. Moore, Minimal two spheres and the topology
of manifolds with positive curvature on totally isotropic two-planes, Annals
of Math., {\bf 127} (1988), 199-227.

\item{11.} S. Nishikawa, Deformation of Riemannian metrics and manifolds
with bounded curvature ratios, Proc. Symp. Pure Math. {\bf 44} (1986),
343--352.

\item{12.} S. Poon, Compact, self-dual manifolds of positive scalar curvature,
J. Diff. Geom., {\bf 24} (1986), 97--132.

\item{13.} S. T. Yau,On the urvature of compact Hermitian manifolds, Inv.
Math.,{\bf25} (1974), 231--239.

\item{14. }J.Cao,Certain 4-manifolds with non-negative curvature,preprint.

\item{16. }R.Hamilton, Four-Manifolds with positive curvature
operator,J.Diff.Geom, {\bf 24}, (1986), 153--179.

\pmf
Institute of Mathematics
\pmf
Hebrew  University
\pmf
Givat Ram 91904 Jerusalem
\pmf
Israel
\pmf
email: simplex\@sunrise.huji.ac.il

\bye